\documentclass[11pt]{article}

\usepackage{amsmath}
\usepackage{amsthm}
\usepackage{amssymb}
\usepackage{graphics}

\title{Geometric renormalisation and Hausdorff dimension for 
loop-approximable geodesics escaping to infinity 
}
\author{
{\em Kurt Falk and Bernd O. Stratmann}\\
\footnotesize{{\sl Universit\"at Bremen, FB 3 - Mathematik,}}\\
\footnotesize{{\sl Bibliothekstra{\ss}e 1, 28359 Bremen, Germany}}\\
\footnotesize{{\tt khf@math.uni-bremen.de, bos@math.uni-bremen.de}}\\
}

\date{}

\def\rz{\mathbb{R}} 
\def\R{\mathbb{R}} 
\def\nz{\mathbb{N}} \def\gz{\mathbb{Z}}
\def\S{\mathbb{S}}
\def\Z{\mathbb{Z}}
 
\def\D{\mathbb{D}}
\def\T{\mathcal{T}}

\def\N{{\cal M}_\H}
\def\H{N}
\def\C{\mathcal{C}}

\def\diam{\mathrm{diam}}

\newtheorem*{theorem*}{Main Theorem}

\newtheorem*{remark*}{Remark}

\begin{document}

\maketitle

\begin{abstract}
\noindent
The main result of this paper is to show that if  $\H$ is a normal subgroup of a Kleinian group $G$ such that $G/\H$ contains a coset which is 
represented by some loxodromic element, then the Hausdorff dimension 
of the transient limit set of $\H$ coincides 
with the Hausdorff dimension of the limit set of $G$. 
This observation extends previous results by Fern\'andez and Meli\'an 
for Riemann surfaces. 
\end{abstract}



\section{Introduction and statement of results}

In this paper we  study  fractal geometric aspects of the
limit set $L(\H)$ of a normal subgroup $\H$ of some
given non-elementary Kleinian group
$G$ acting on $(m+1)$-dimensional hyperbolic space $\D^{m+1}$.
We always assume that $G/\H$ contains
a coset which is  represented by some loxodromic isometry $\gamma \in G$.  
It is well known that  in this situation $L(\H)$ coincides with
the limit set $L(G)$ of the larger group $G$. 
However, a comparison  of finer aspects of these two limit sets  
usually turns out to be far more involved, as can be seen, for 
instance, in the work of Brooks \cite{brooks}
and Rees \cite{rees1, rees2}.\\
In this paper we investigate the set $L_{t}(\H)$ of  directions 
at some arbitrary point $z$ on the
manifold ${\cal M}_{\H}$ associated with $\H$ for which  the 
resulting geodesic movement on ${\cal M}_{\H}$ 
eventually escapes from every compact region on ${\cal M}_{\H}$, but 
which is nevertheless  contained in the 
$\epsilon$-neighbourhood of some sequence of closed loops starting 
and ending at
$z$ on ${\cal M}_{\H}$, for each $\epsilon>0$. That is, 
we consider the transient limit set $L_{t}(\H)$ of $\H$, given by
$$
L_{t}(\H):=\{\xi \in L(\H): \lim_{r\to \infty} d(\xi(r), \H(0)) 
=\infty\}.
$$ 
Here, $d$ refers to the hyperbolic metric in $\D^{m+1}$ 
and $\xi(r)$ denotes a $\R_{+}$-parametrisation of the  
geodesic ray from the origin  to $\xi$. \\  The 
following theorem gives the main result of this paper.

\begin{theorem*}
Let $G$ and $\H$  be given as above. We then  have 
\[
\dim_H(L_t(\H)) = \dim_H(L(G)).
\]
\end{theorem*}

This theorem gives a extension 
of results by Fern\'andez and Meli\'an \cite{fermel01}, 
who studied the set ${\cal E}$ of
escaping, not necessarily loop-approximable directions on a
complete oriented non-compact Riemann surface ${\cal R}$ with 
fundamental group $\Gamma$. 
That is, 
${\cal E}:= \{\xi \in  \S^{1}: \lim_{r \to \infty} \rho 
 (\xi(r), \Gamma(0))  = \infty\}$, 
where $\S^{1}$ refers to the boundary at infinity of $\D^{2}$
and $\rho$ denotes the hyperbolic metric in $\D^{2}$. 
The main result of \cite[Theorem~1]{fermel01} was to establish the 
following tricothomy.
 \begin{itemize}
     \item[  (i)] If ${\cal R}$ has finite area, then ${\cal E}$ 
     is countable.
     \item[ (ii)] If Brownian motion on ${\cal R}$ is transient, then
     ${\cal E}$ 
     has full Lebesgue measure.
     \item[(iii)] If ${\cal R}$ has infinite area and Brownian motion 
     is recurrent, then ${\cal E}$  
     has zero Lebesgue measure, but its Hausdorff dimension is equal 
     to $1$.
     \end{itemize}
Therefore, for hyperbolic manifolds which are normal coverings of 
some hyperbolic manifold and which posess a loxodromic representative $\gamma 
\in G/ \H$, 
our Main Theorem extends the results by Fern\'andez and Meli\'an to
arbitrary dimensions and to the situation where the boundary at 
infinity of hyperbolic space is replaced by the limit set of the 
fundamental group. However, let us emphasize 
that our proof does require the existence of a loxodromic $\gamma 
\in G/ \H$,  and hence our 
extension is restricted to normal coverings with this property.

\begin{figure}
\caption{Dynamics in $\N$.}
$\, $
{\par
    \centering
    \includegraphics{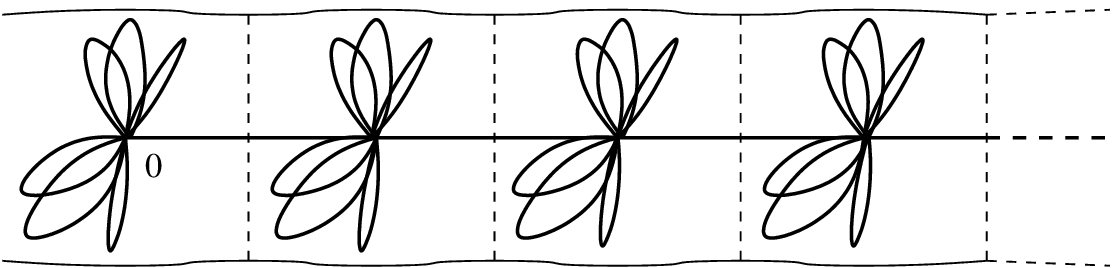}
    \label{mnf}
    \par}
\end{figure}

Our proof hinges on two renormalisation procedures.  
These are used to locate a certain family of  subsystems within 
the loop-approximable, non-recurrent dynamics on the manifold $\N$.
At the boundary of  the universal covering space these subsystems are described by a 
family of Cantor sets  contained in $L(\H)$, and here the key 
observation is that this 
family contains Cantor sets whose Hausdorff dimension is 
arbitrarily close to the Hausdorff dimension of $L(G)$.
These Cantor sets are constructed inductively using the
following two renormalisation procedures. The first of these employs a 
well-known construction by  Bishop and Jones  for computing
the Hausdorff dimension of bounded dynamics (see \cite{bijo97},
\cite{str01a}). 
This construction gives rise to a certain weighted scaling law, which 
we then apply  a sufficient number of times in order to prepare for the second
renormalisation step. 
On the manifold $\N$ this first step of the overall construction corresponds
to a family of well separated quasi-geodesics   within a bounded 
region of $\N$,
each starting at the same point.
The second renormalisation procedure consists  in  prolonging  
each  of these quasi-geodesics  by a long geodesic segment which is
 contained in the projection of the axis $A_\gamma$ of $\gamma$ onto 
 $\N$, and which is chosen  such that 
it leads out of the bounded region which contained the original  
quasi-geodesics.  
The resulting dynamical behaviour on $\N$ is sketched in
Figure~\ref{mnf}. 

We assume that the reader is familiar with the proof of  
Bishop and Jones'  result on the relationship between the 
exponent of convergence $\delta(G)$ of a non-elementary Kleinian 
group $G$ and the Hausdorff dimension of the radial limit set 
$L_r(G)$ of $G$ (see \cite{bijo97} and for a more detailed proof 
\cite{str01a}).
 Here, the reader might like to recall that   the radial limit set 
$L_r(G)$ represents those limit points $\xi$ for                
   which the projection of a hyperbolic ray towards $\xi$ returns infinitely   
   often to some compact part of the manifold associated to $G$.    

\section{Preliminaries}
Throughout, let $G$ be a non-elementary Kleinian group acting on 
$(m+1)$-dimensional hyperbolic space $\D^{m+1}$. Also, let $\H$ 
be a normal subgroup of $G$ such that $G/\H$ contains a coset 
$[\gamma]$, for some loxodromic $\gamma \in G$.
Moreover, we always assume without loss of generality  that
$0\in \D^{m+1}$ is an element of the axis $ A_{\gamma}$
of $\gamma$, and we let $\eta_-$ ($\eta_+$ resp.) refer 
to the repulsive (attractive resp.) fixed point of $\gamma$.
 Next,  recall that  to any arbitrary  Kleinian group $\Gamma$ we can associate its 
truncated Poincar\'e series ${\cal P}_{t}(\Gamma,s,w) $, as well as
its Poincar\'e series ${\cal P}(\Gamma,s,w)$. These series are
given, for $s,t \in \R_{+}$ and $w \in \D^{m+1}$, by 
\[ 
{\cal P}_{t}(\Gamma,s,w) := \sum_{h \in \Gamma \atop d(w,h(w)) \leq t} 
e^{-sd(w,h(w))}, \, \, \hbox{ and } \, \, 
{\cal P}(\Gamma ,s,w):= \lim_{t \to 
\infty} {\cal P}_{t}(\Gamma,s,w) .
\]
The abzissa of convergence of the infinite series ${\cal P}(\Gamma,s,w)$
is called the exponent of convergence of $\Gamma$, and  it will be denoted 
by $\delta(\Gamma)$. Here, note that Bishop and Jones \cite{bijo97}
(see also \cite{str01a})
showed that if $\Gamma$ is non-elementary, then we always have that
\[ 
\delta(\Gamma)= \dim_{H}(L(\Gamma)).
\]
Also, we will make use of the following standard facts and notations for the
Poincar\'e model $\left(\D^{m+1},d\right)$ of the 
$(m+1)$-dimensional hyperbolic space. 
Let $B(w,r)$ refer to the hyperbolic ball centred
at $w$ of radius $r$,  and let $\Pi:\D^{m+1} \to  \S^{m}$ 
denote the radial projection from the origin
to the boundary $\S^{m}$ of hyperbolic space. That is,
for $E \subset \D^{m+1}$ we have $\Pi(E):=\{\xi \in \S^{m}: s_{\xi}
\cap E \neq \emptyset\}$, where $s_{\xi}$ refers to the Euclidean 
straight line between
the origin and $\xi$.
Also, we require the following analogue of Pythagoras's Theorem for 
hyperbolic triangles. For this, consider a hyperbolic triangle
with sides of finite lengths $a$, $b$ and $c$, and  with 
$\alpha_{0}\in (0,\pi)$ denoting the angle opposite to the side 
of length $a$.
A straightforward application of the hyperbolic cosine rule (see e.g.
\cite{be}) then gives that 
there exists a constant $K>0$, depending only on $\alpha_{0}$, such that 
$$
b+c-K \leq a \leq b+c.
$$
Moreover, we require the following additional  observation from elementary 
hyperbolic geometry. For this, consider some arbitrary hyperbolic 
geodesic $A \subset \D^{m+1}$ which does not contain the origin,
and let $s_{\eta}$ denote the geodesic ray connecting the
origin with one of the endpoints $\eta \in \S^{m}$  of $A$. Also, let 
$\hat{z}_{A}$ refer to the summit of the geodesic $A$. That is,
$\hat{z}_{A}$ is uniquely determined by 
$d(0, \hat{z}_{A}) = \min \{ d(0,w): w \in A\}$. 
A straightforward exercise in hyperbolic geometry then shows that 
there exists a universal constant $\tau>0$ such that
\begin{eqnarray}
\label{axis} 
\min\{ d( w,\hat{z}_{A}):  w \in s_{\eta} \} < \tau.
\end{eqnarray}
In fact,  an elementary calculation shows
that $\tau$ is equal to $\log(1+\sqrt{2})$, which is often referred to 
as Schweikart's constant.
 
For a sequence $(w_n)_{n \in \nz_0}$ of distinct points 
in $\D^{m+1}$, let
$$
[w_0, w_1, w_2	, \ldots ]
$$
denote the quasi-geodesic path obtained by connecting  $w_n$ and 
$w_{n+1}$ with the unique geodesic arc between them, 
for each $n \in \nz$. A standard observation from 
hyperbolic geometry then shows that if the lengths of these geodesic segments are 
uniformly bounded away from $0$, and if  each of the angles between  adjacent 
geodesic segments 
is uniformly bounded from below by some $\alpha_{0}>0$, then 
$[w_0, w_1, w_2	, \ldots ]$ is  a quasi-geodesic
ray towards a unique point at infinity. That is, each $w_n$ is, 
with respect to the hyperbolic metric, uniformly 
bounded (depending on $\alpha_{0}$)
 away from the geodesic ray from $w_0$ towards the 
uniquely determined limit  at the boundary at  infinity of the 
sequence $(w_n)$.
  
Finally, note that we use the common notation $a_n \asymp b_n$
if  two sequences of positive real numbers $a_n$ and $b_n$ are 
comparable, that is,  if the ratio $a_n/b_n$ is uniformly
bounded from  below by $1/c$ and from above by $c$, for some $c>1$ 
and for all $n \in \nz$.

\section{The two renormalisation procedures}\label{3}

Let us begin with by giving our first renormalisation procedure.
Here, a {\em geodesic $\H$-tree ${\cal T}(z)$ rooted at $z$} refers 
to an infinite tree whose set of vertices $V(\T(z))$ is contained 
in $\H(z)$ and whose edges are finite geodesic segments 
between the vertices, such that each vertex $u\in {\cal T}(z)$ has a
finite set ${\cal S}(u)$ of successors of cardinality at least $2$  and  
such that each element in  $V(\T(z))\setminus \{z\}$ has a unique
predecessor. 

In the following,  let $z_{n}$ be defined by $z_{n}:=\gamma^{n}(0)$, for  each 
$n \in \gz$.
Note that our first renormalisation procedure
is well known for the special case in which each $z_{n}$ lies
in the orbit $\H(0)$. In this situation, its 
outcome  has already been obtained in 
\cite[Proposition 3.5]{str01a}. The novelty here is that for the normal 
subgroup $\H$ the result of \cite[Proposition 3.5]{str01a} continues to 
hold for each element of the orbit $\{z_{n}: n \in \Z\}$ 
of the origin under $\langle \gamma\rangle$.

\vspace{2mm}

{\bf Recurrent Renormalisation Procedure ([RRP]).} $\,$

{\em
For each $0< s < \delta(\H)$,  there exist
$\kappa>0$, $\ell_{s}>0$ and $K_{s}>1$ such that for each $h \in \H$ 
and  $n \in \Z$ there exists a geodesic $\H$-tree 
$\T=\T_{s}(h_{0}(z_n))$ rooted at $h(z_n)$ 
with the following properties.
\begin{itemize}
\item[  (i)] If $u \in V(\T)$, then $\Pi(B(v, \kappa)) \subset 
\Pi(B(u, \kappa))$ for each $v \in {\cal S}(u)$. 
\item[ (ii)] If $v \in {\cal S}(u)$ for some $u \in V(\T)$, then $d(u,v) \leq 
\ell_{s}$.
\item[(iii)] If $v,w \in {\cal S}(u)$ for some $u \in V(\T)$, then
$\exp(d(0,v)) \asymp \exp(d(0,w))$ and 
$\Pi(B(v,\kappa)) \cap \Pi(B(w,\kappa)) = \emptyset$. 
\item[ (iv)]For each $u \in V(\T)$ we have
$$
\sum_{v \in {\cal S}(u)} (\diam(\Pi(B(v,\kappa))))^s \geq 
K_{s} \, (\diam(\Pi(B(u,\kappa))))^s .
$$
\end{itemize}
We will say that the  so derived family 
$\{B(v,\kappa):  v \in {\cal S}(u)\}$ is 
obtained by applying the recurrent renormalisation procedure 
to $u \in V(\T)$. 
}

\begin{proof}
As already mentioned before, for $n=0$ the assertion in this procedure
has been obtained in \cite[Proposition 3.5]{str01a} (see also 
\cite{bijo97}), and we refer to these papers for the proof in this case.
In fact, note that the main idea of the proof  
of \cite[Proposition~3.5]{str01a}  consists of a geometrization 
of the rate of increase, for $t$ tending to infinity,  of the truncated Poincar\'e series ${\cal P}_{t}(H,s,0)$ for $s < \delta(H)$.
For the proof of the general situation, that is, 
for some arbitrary $n \in \Z$,
note that the value of the truncated Poincar\'e series associated with 
$\H$ does not change if we exchange the observation point $z_{0}=0$ by
some arbitrary point  in $\{z_{n}:n \in \Z\}$. 
More precisely,  since $\H$ is normal in $G$,  
we have for each  $n \in \gz$ and $s,t \in \rz_{+}$,
\begin{eqnarray*}
{\cal P}_{t}(\H,s, z_{n}) 
& = & 
\sum_{h \in \H \atop d(z_{n}, h (z_{n})) \leq t} 
e^{-sd(z_{n}, h(z_{n}))} 
  =
\sum_{h \in \H \atop d(\gamma^{n}(0), h \gamma^{n}(0)) \leq t}
e^{-sd(\gamma^{n}(0), h \gamma^{n}(0))} \\
& = &
\sum_{\gamma^{-n} h \gamma^{n} \in \H \atop 
d(0,\gamma^{-n} h \gamma^{n}(0)) \leq t}
e^{-sd(0,\gamma^{-n} h \gamma^{n}(0))} 
  =  
\sum_{h \in \H \atop d(0, h (0)) \leq t} e^{-sd(0, h (0))} \\
& = &
{\cal P}_{t}(\H,s, 0).
\end{eqnarray*}
Using this observation, the assertion now follows from a straightforward 
adaptation of the arguments in the proof of 
\cite[Proposition 3.5]{str01a}. 
\end{proof}

\vspace{2mm}

{\bf Transient Renormalisation Procedure ([TRP]).} $\,$

{\em 
For $0< s< \delta(\H)$,  $n \in \Z$ and $h_{0}\in \H$, let
$\T=\T_{s}(h_{0}(z_n))$ denote the geodesic $\H$-tree obtained in the recurrent
renormalisation procedure {\em [RRP]}.  
Then there exists a  constant $0< k_{\gamma} < 1$  such 
that for each $q \in \nz_{0}$ sufficiently large and for each  $h \in \H$ with
$h(z_n) \in V(\T) \setminus \{ h_0(z_n) \}$,  we have that
$$
\diam(\Pi(B(h(z_{n+q}),\tau)) )\geq k_{\gamma}^{\: q} \,
\diam(\Pi(B(h(z_n),\tau))),
$$
where $\tau:=\log(1 + \sqrt{2})$. We will say that the  ball 
$B(h(z_{n+q}),\tau)$ is obtained by starting at $h(z_n)$ and 
applying the transient renormalisation procedure $q$ times. }

\begin{proof}

\begin{figure}
\caption{The location of $h(\eta_+)$.}
$\, $
{\par
    \centering
    \includegraphics{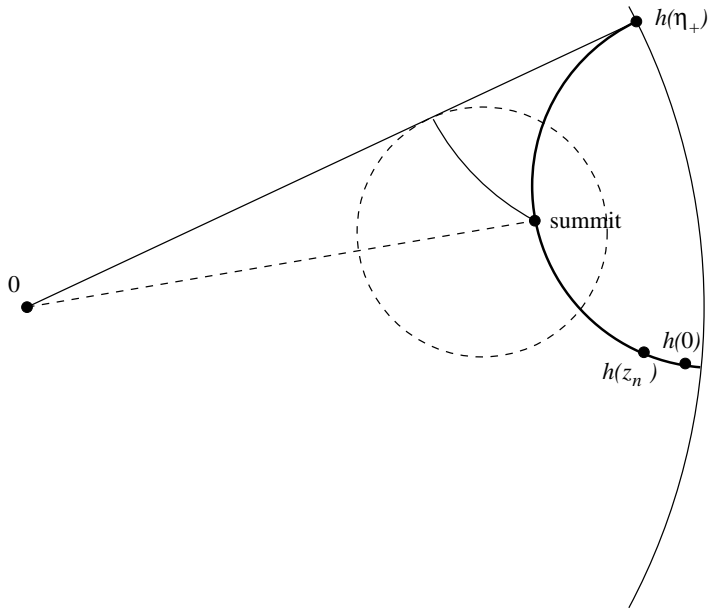}
    \label{etaplus}
    \par}
\end{figure}

Let $\T=\T_{s}(h_{0}(z_n))$ and $h \in \H$ be given as stated 
in the renormalisation procedure. Let us first show that $h(z_{n})$
lies always close to the summit $\hat{z}_{h(A_{\gamma})}$ of the 
geodesic $h(A_{\gamma})$,  
the image of the axis $A_{\gamma}$ under $h$. Indeed, 
since $h(z_n) \in V(\T) \setminus 
\{ h_{0}(z_n) \}$,  the statement in (1)  of [RRP] implies that there exists 
$u\in V(\T)$ such that  
$\Pi(B(h(z_{n}),\kappa)) \subset \Pi(B(u,\kappa))$ (see also Figure~\ref{etaplus}).

The elementary observation in (\ref{axis}) shows that the distance from
the summit of a geodesic to each of the two rays from the origin to the
endpoints of the geodesic is 
less than  $\tau$. Thus, by construction, 
we have that
$$
h(\eta_+) \in \Pi(B(\hat{z}_{h(A_{\gamma})},\tau)).
$$

Recall that in the recurrent renormalisation procedure [RRP] we have already 
derived the existence of the parameter $\ell_{s}$, which is the upper 
bound of the lengths of the edges in the tree $\T$.

Next, consider the geodesic  $h(A_\gamma)$ containing the points
$h(0)$ and $h(z_n)$; one of its endpoints will be $h(\eta_+)$.
Assume, by way of contradiction, that the distance between $h(z_n)$
and the summit $\hat{z}_{h(A_{\gamma})}$ of $h(A_\gamma)$ is larger than 
$2 \ell_s + \tau$.
Projecting onto the manifold $\N$ associated to $\H$ and using the
hyperbolic triangle inequality, we obtain a 
contradiction to the fact that 
$h(z_n) \in V(\T) \setminus \{ h_{0}(z_n) \}$.
It  immediately follows that 
$$
\diam(\Pi(B(\hat{z}_{h(A_{\gamma})},\tau))) \asymp 
\diam(\Pi(B(h(z_n),\tau))),
$$
where the comparability  constant depends only  on the distance 
$2 \ell_s + \tau$, and therefore,
only on $\H$ and $s$.
The statement now follows by applying $h\gamma^{q}h^{-1}$ to the ball $B(h(z_n),\tau)$, 
which immediately gives that
\begin{eqnarray*}
\diam(\Pi(B(h(z_{n+q}),\tau)) ) 
&\geq &
k_{\gamma}^{\: q} \,
\diam(\Pi(B(h(z_n),\tau))),
\end{eqnarray*}
where $k_{\gamma} \asymp \exp(-d(0, \gamma(0)))$. 
\end{proof}

Let us remark that the constant $\kappa$ in [RRP] and the constant
$\tau$ in [TRP]  are
independent of each other.  Also, the statements in [RRP] 
continue to hold if we replace $\kappa$ by a smaller positive 
number, and the same holds for $\tau$ in [TRP].
Therefore, for the remainder of this paper, when applying [RRP] and 
[TRP],  we  use 
$$\sigma:= \min \{\kappa, \tau \}$$
instead of $\kappa$ and $\tau$.

\section{Proof of the theorem}
Let us first observe that it is sufficient to prove the assertion in the 
theorem for the case in which
\[ 
\delta(\H) = \dim_{H}(L(G)).
\]
Indeed, this can immediately be seen by way of contradiction as 
follows. Suppose that $\delta(\H) < \dim_{H}(L(G))$.  Since diminishing 
a set by a subset of smaller Hausdorff dimension does not alter the 
Hausdorff dimension of that set and  using the well known fact that
$\delta(\H) = \dim_{H}(L_{r}(\H))$  
(see \cite{bijo97} and \cite{str01a}),  we obtain
\[
\dim_{H}(L(G)) = \dim_{H}(L(\H)) = \dim_{H}(L_{t}(\H) \cup L_{r}(\H)) =
\dim_{H}(L_{t}(\H)).
\]
Therefore, we can now assume, without loss of generality, that
$$
\delta(\H) = \dim_{H}(L(G))= \delta(G).
$$
The rough strategy for proving the main theorem in this case is
as follows. For  some arbitrary given $0<s<\delta(\H)$, we construct
a certain Cantor set $\C_s \subset L_t(\H)$, and then show that
$\dim_H(\C_s) \geq s$. By the arbitrary choice of $s$, the theorem then
follows.
The idea of  the Cantor set construction is to  start at the origin and then to perform an alternating inductive process 
using both renormalisation procedures. The building block of this process is that we first apply the recurrent renormalisation procedure [RRP] sufficiently many times 
until the resulting power of $K_{s}$ is large enough 
(in fact, this number of times depends on the outcome of the step to come). 
After that, we perform the transient  renormalisation procedure [TRP] 
sufficiently many times, without loosing the control on the distortion
(in particular, this step will guarantee that out Cantor set contains only transient limit points). 
More precisely, 
let $0<s< \delta(\H)$ be given. Then $s$ determines the width 
$\ell_{s}$ of the recurrent renormalisation procedure [RRP]. 
Having fixed $\ell_{s}$, we choose  $q \in \nz$, the  number
of times we are going to
apply the transient  renormalisation procedure [TRP], so that
\begin{eqnarray}
\label{transient}
q \, d(0,\gamma(0))  \geq 4 \ell_{s}.
\end{eqnarray}
This choice of $q$ will guarantee that the Cantor set $\C_s$ we are going to construct  will  be contained in $L_t(\H)$.   
Finally, we choose $p \in \nz$ to  be minimal with respect to the property 
\begin{eqnarray}
\label{constants}
K_{s}^{\,p} \, k_{\gamma}^{\,q} > 1.
\end{eqnarray}
Let us now come to the explicit construction of $\C_s$. As already mentioned, the construction starts at the origin, and 
we set $T_ {0}(z_{0})= T_{0}(0):= \{0\}$. Then, the first step is to apply the recurrent renormalisation procedure [RRP] $p$ times, starting at
 $z_{0}$. 
According to [RRP], this gives rise to a set of hyperbolic balls 
$B(v,\sigma)$ whose radial projections to the boundary $\S^{m}$ are pairwise disjoint and of comparable diameter. 
The set of centres of these balls in $\D^{m+1}$ will be denoted 
by $R_ {p}(z_{0})$.  Then, the second step is to apply  the   transient  renormalisation 
procedure [TRP] $q$ times to each element in $R_ {p}(z_{0})$.  The set of centres of the
so obtained hyperbolic balls will be denoted by 
$T_ {p}(z_{q})$.  This represents the start of the induction, and we then continue
 as follows. 
 Assume that the sets $R_ {np}(z_{(n-1)q})$ and $T_ {np}(z_{nq})$ have been constructed,  for some $n \in \nz$.
 To each of the points in $T_ {np}(z_{nq})$ we then  
apply  the recurrent renormalisation procedure [RRP] $p$ times. The set of centres of these so obtained hyperbolic balls 
gives the set $R_ {(n+1)p}(z_{nq})$. Next, we apply the   transient  renormali\-sation 
procedure [TRP] $q$ times to each of the elements in $R_ {(n+1)p}(z_{nq})$.  The set of centres of these
so obtained hyperbolic balls gives rise to the set  
$T_ {(n+1)p}(z_{(n+1)q})$. 
\begin{figure}[ht]
\caption{The inductive construction of $\C_s$.}
$\, $
{\par
    \centering
    \includegraphics{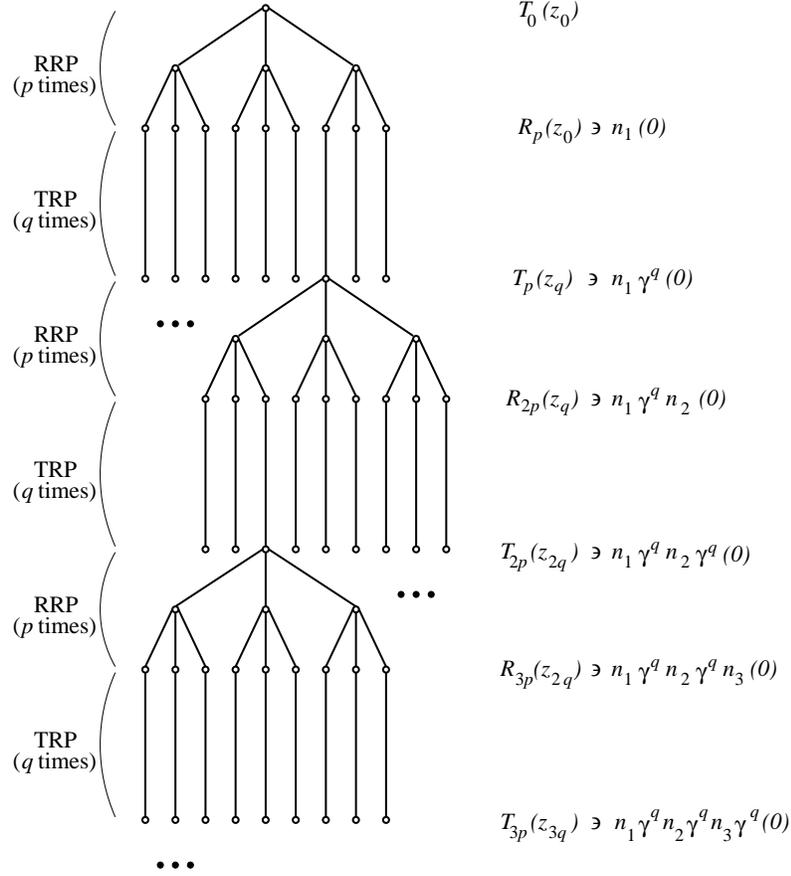}
    \label{induct}
    \par}
\end{figure}
This finishes our alternating inductive argument 
(see also Figure~\ref{induct}), and we can now use it to define our
desired Cantor set $\C_s$ by
$$
\C_s := \bigcap_{n \in \nz} \; \bigcup_{v\in T_ {np}(z_{nq})} 
\Pi(B(v,\sigma)).
$$
Here,  $\sigma>0$ refers to the constant which we specified at the end of Section \ref{3}. \\
Next, observe that in this Cantor set construction we have  good control  over the distortion, when going from one generation 
in the construction to the next. That is,  by using
[RRP] $(iv)$,  [TRP] and the condition in (\ref{constants}), we have the following crucial estimate, for each $n \in \nz$,
\begin{eqnarray*}
\sum_{v\in T_ {np}(z_{nq})}  (\diam(\Pi(B(v,\sigma))))^s  
& \geq & 
k_{\gamma}^{q} \sum_{v\in R_ {np}(z_{(n-1)q}) } 
(\diam(\Pi(B(v,\sigma))))^s \\ 
& \geq & 
K_{s}^{\,p} \, k_{\gamma}^{q} \; \sum_{v\in T_ {(n-1)p}(z_{(n-1)q}) } 
(\diam(\Pi(B(v,\sigma))))^s\\
& > &
\sum_{v\in T_ {(n-1)p}(z_{(n-1)q}) } 
(\diam(\Pi(B(v,\sigma))))^s .
\end{eqnarray*}
Using a straightforward generalisation of the folklore arguments from fractal geo\-metry of \cite[Lemma 2.5]{str01a} and  
\cite[Corollary 2.6]{str01a}, the latter estimate  
immediately gives  that
$$
\dim_H(\C_s) \geq s.
$$
It remains to show that the set $\C_s$ is contained in $L_{t}(\H)$. 
For this, note that, by viewing 
the construction of $\C_s$ from within $\D^{m+1}$, the set  $\C_s$ gives rise to a geodesic $G$-tree which is rooted at the origin
and whose vertex set is  equal to 
$\bigcup_{n \in \nz} \left(T_ {(n-1)p}(z_{(n-1)q}) \cup 
R_ {np}(z_{(n-1)q}) \right)$.  
By construction, this tree has the property that the lengths of the constituting geodesic edges and the angles formed by adjacent edges
are uniformly bounded away from zero. Therefore, each path in this tree starting at the origin  is a quasi-geodesic heading towards a uniquely determined point at infinity.  Clearly, the projection of each of these quasi-geodesics  onto the manifold ${\cal M}_{\H}$ gives
some piecewise geodesic movement in ${\cal M}_{\H}$ 
which has the following properties. If an edge in the tree starts at 
a vertex in $T_ {(n-1)p}(z_{(n-1)q})$ and
ends at a point in $R_ {np}(z_{(n-1)q})$, then in ${\cal M}_{\H}$ 
this edge is represented by a geodesic loop of hyperbolic length 
at most $\ell_s$. Obviously, this loop must then be contained in 
a bounded region of  ${\cal M}_{\H}$ of diameter at most $\ell_s$. 
Whereas, if an edge in the tree starts at a vertex in 
$R_ {np}(z_{(n-1)q})$ and
ends at a point in $T_ {np}(z_{nq})$, using [TRP],  we then have that  
in ${\cal M}_{\H}$ this edge represents a geodesic segment in 
${\cal M}_{\H}$ which starts in the previous bounded region and then  
heads straight towards the end of ${\cal M}_{\H}$ associated with 
the attractive fixed point $\eta_+$ of $\gamma$. 
Moreover, the condition in (\ref{transient}) guarantees that
the hyperbolic length of that segment is at least equal to $4 \ell_s$, 
and this  shows that it's end point  is separated by at least $3 \ell_s$
from the previous bounded region. 
\hfill $\Box$

\end{document}